\input amstex
\documentstyle{amsppt}
\NoBlackBoxes
\magnification 1200
\hsize 6.5truein
\vsize 8truein
\topmatter
\title A Slight Improvement to Garaev's Sum Product estimate \endtitle
\author  Nets Hawk Katz and Chun-Yen Shen \endauthor
\affil Indiana University \endaffil
 \thanks 
  The first  author was supported by NSF grant DMS 0432237.
\endthanks 
\endtopmatter

\head \S 0 Introduction \endhead

Let $A$ be a subset of $F_p$, the field of $p$ elements with $p$ prime.

We let
$$A+A=\{a+b: a \in A, b \in A \},$$
and
$$AA=\{ab: a \in A, b \in A \}.$$

It is fun (and useful) to prove lower bounds on $\max(|A+A|,|AA|)$ (see e.g. [BKT],[BGK],[G]).
Recently, Garaev [G] showed that when $|A| < p^{{1 \over 2}}$ one has the estimate
$$\max(|A+A|,|AA|) \gtrapprox |A|^{{15 \over 14}}.$$

By using Plunneke's inequality in a slightly more sophisticated way, we improve this exponent to
${{14 \over 13}}$. We believe that further improvements might be possible through aggressive use
of Ruzsa covering.

\head \S 1 Preliminaries \endhead

Throughout this paper $A$ will denote a fixed set in the field $F_p$ of $p$ elements with $p$
a prime. For $B$, any set, we will denote its cardinality by $|B|$.

Whenever $X$ and $Y$ are quantities we will use
$$X \lesssim Y,$$
to mean 
$$X \leq C Y,$$
where the constant $C$ is universal (i.e. independent of $p$ and $A$). The
constant $C$ may vary from line to line.
We will use
$$X \lessapprox Y,$$
to mean
$$X \leq C (\log |A|)^{\alpha} Y,$$
where $C$ and $\alpha$ may vary from line to line but are universal.

We state some preliminary lemmas, mostly those stated by Garaev, but occasionally with
different emphasis.

The first lemma is a consequence of the work of Glibichuk and Konyagin [GK]

\proclaim{Lemma 1.1} Let $A_1 \subset F_p$ with $1 < |A_1| < p^{{1 \over 2}}$.
Then for any elements $a_1,a_2,b_1,b_2$ so that
$${b_1-b_2 \over a_1-a_2} + 1 \notin {A_1-A_1 \over A_1-A_1},$$
we have that for any $A^{\prime} \subset A_1$ with $|A^{\prime}| \gtrsim |A_1|$
$$|(a_1-a_2) A^{\prime} + (a_1 - a_2) A^{\prime} + (b_1 -b_2) A^{\prime}|
 \gtrsim |A_1|^2.$$
In particular such $a_1,a_2,b_1,b_2$ exist unless ${A_1-A_1 \over A_1-A_1}=F_p$.
In case ${A_1-A_1 \over A_1-A_1}=F_p$, we may find $a_1,a_2,b_1,b_2 \in A_1$
so that
$$|(a_1-a_2) A_1 + (b_1-b_2) A_1| \gtrsim |A_1|^2.$$
\endproclaim

\demo{Sketch of Proof} If ${A_1-A_1 \over A_1-A_1} \neq F_p$, it is immediate that there exist
$a_1,a_2,b_1,b_2 \in A_1$ with $1+{b_1-b_2 \over a_1-a_2} \notin {A_1 - A_1 \over A_1-A_1}$.
This automatically implies
$$|(a_1-a_2) A^{\prime} + (a_1 - a_2) A^{\prime} + (b_1 -b_2) A^{\prime}|
 \gtrsim |A_1|^2.$$ (See [GK]. If $x \notin {A_1-A_1 \over A_1 - A_1}$ then each element of $A_1+xA_1$
has but one representative $a+xa^{\prime}$.
On the other hand if
$${A-A \over A-A} = F_p,$$
then one can find $a_1,a_2,b_1,b_2 \in A_1$ so that ${a_1-a_2 \over b_1 -b_2}$ has at most $|A|^2$
representatives as ${a_3-a_4 \over b_3-b_4}$ with $a_3,a_4,b_3,b_4 \in A$
 which implies that $|A+{a_1-a_2 \over b_1 -b_2}A|$ is large. Again, for more details see [GK].
\qed \enddemo

The following two lemmas, quoted by Garaev, are due to Ruzsa, may be found in [TV].
The first is usually referred to as Rusza's triangle inequality. The second is
a form of Plunneke's inequality.

\proclaim{Lemma 1.2} For any subsets $X,Y,Z$ of $F_p$ we have
$$|Y-Z| \leq {|Y-X||X -Z| \over |X|}.$$ \endproclaim

\proclaim{Lemma 1.3} Let $X,B_1,\dots,B_k$ be any subsets of $F_p$ with
$$|X+B_i| \leq \alpha_i |X|,$$
for $i$ ranging from 1 to $k$. Then there exists $X_1 \subset X$
with
$$|X_1+B_1 + \dots + B_k| \leq \alpha_1 \dots \alpha_k |X_1|. \tag (1.1)$$
\endproclaim

We record a number of Corollaries. The first two can be found in [TV]. The last
one, we first became aware of in the paper of Garaev.

\proclaim{Corollary 1.4} Let $X,B_1,\dots,B_k$ be any subsets of $F_p$.
Then
$$|B_1 + \dots + B_k| \leq {|X+B_1| \dots |X+B_k| \over |X|^{k-1}}.$$
\endproclaim

\demo{Proof} Simply bound $|B_1 + \dots + B_k|$ by $|X_1+B_1 + \dots + B_k|$
and $|X_1|$ by $|X|$. \qed \enddemo

Corollary 1.4 is somewhat wasteful in that $X_1$ is unlikely to be both a
singleton element and a set with the same cardinality as $X$. By applying
Lemma 1.3 iteratively.

\proclaim{Corollary 1.5} Let $X,B_1,\dots,B_k$ be any subsets of $F_p$.Then
there is $X^{\prime} \subset X$ with $|X^{\prime}| > {1 \over 2} |X|$ so that
$$|X^{\prime}+ B_1 + \dots B_k| \lesssim {|X+B_1| \dots |X+B_k| \over |X|^{k-1}}.$$
\endproclaim

\demo{Proof} Observe that for any $Y \subset X$ with $|Y| \geq {|X| \over 2}$,
we have that
$${|Y+B_i| \over |Y|} \lesssim {|X+B_i| \over |X|}.$$

Now recursively apply Lemma 1.3. That is, first apply it to $X,B_1,\dots,B_k$
obtaining a set $X_1$ satisfying
$$|X_1+B_1+\dots+B_k| \lesssim {|X+B_1| \dots |X+B_k| \over |X|^{k}}|X_1|.$$
If $|X_1| > {1 \over 2} |X|$ then stop and let $X^{\prime}=X_1$. Otherwise
apply Lemma 1.3 to $X \backslash X_1, B_1 ,\dots,B_k$. Proceeding recursively
if $|X_1 \cup \dots \cup X_{j-1}| \geq {1 \over 2} |X|$ then set
$$X^{\prime} = X_1 \cup \dots \cup X_{j-1},$$
otherwise obtain the inequality
$$|X_j+B_1+\dots+B_k| \lesssim {|X+B_1| \dots |X+B_k| \over |X|^{k}}|X_j|.$$
Summing all the inequalities we obtained before stopping gives us the desired
result.
\qed \enddemo

\proclaim{Corollary 1.6} Let $A \subset F_p$ and let $a,b \in A$. Then we have
the inequalities
$$|aA+bA| \leq {|A+A|^2 \over |aA \cap bA|},$$
and
$$|aA-bA| \leq {|A+A|^2 \over |aA \cap bA|}.$$
\endproclaim

\demo{Proof}
To get the first inequality, apply Corollary 1.4 with $k=2$, $B_1=aA$, $B_2=bA$, and
$X=aA \cap bA$.

To get the second inequality, apply Lemma 1.2 with $Y=aA$, $Z=-bA$ and $X=-(aA \cap bA)$.
\qed \enddemo

\head \S 2 Modified Garaev's inequality  \endhead

In this section, we slightly modify Garaev's argument to obtain

\proclaim{Theorem 2.1} Let $A \subset F_p$ with $|A| < p^{{1 \over 2}}$ then
$$\max(|AA|, |A+A|) \gtrapprox |A|^{{14 \over 13}}.$$ \endproclaim

\demo{Proof}

Following Garaev, we observe that
$$\sum_{a \in A} \sum_{b \in A} |aA \cap bA| \geq {|A|^4 \over |AA|}.$$
Therefore, we can find an element $b_0 \in A$, a subset $A_1 \subset A$
and a number $N$ satisfying 
$$|b_0 A \cap a A| \approx N,$$
for every $a \in A_1$. Further
$$N \gtrapprox {|A|^2 \over |AA|}, \tag 2.1$$
and 
$$|A_1| N \gtrapprox {|A|^3 \over |AA|}. \tag 2.2 $$

Now there are two cases. In the first case, we have
$${A_1-A_1 \over A_1-A_1} = F_p.$$
If so, applying Lemma 1.1,
we can find $a_1,a_2,b_1,b_2 \in A_1$ so that
$$|A_1|^2 \lesssim |(a_1-a_2) A_1 + (b_1-b_2) A_1| \leq |a_1 A - a_2 A + b_1 A - b_2 A|.$$

Applying Corollary 1.4 with $k=4$ and with $B_1= a_1 A$, with $B_2=-a_2A$ with
$B_3 = b_1 A$, with $B_4=-b_2 A$, and with $X=b_0 A$. Then we apply
Corollary 1.6 to bound above $|X+B_j|$. This yields

$$|A_1|^2 \lesssim {|A+A|^8 \over N^4 |A|^3},$$
or
$$|A_1|^2 N^4 |A|^3 \lessapprox |A+A|^8 .$$
Applying (2.2), we get
$$N^2 |A|^9 \lessapprox |A+A|^8 |AA|^2. \tag 2.3$$
and applying (2.1), we get
$$|A|^{13} \lessapprox |A+A|^8 |AA|^4. \tag 2.4$$
The estimate (2.4) implies that 
$$\max(|A+A|,|AA|) \gtrapprox |A|^{{13 \over 12}}
\gtrapprox |A|^{{14 \over 13}},$$
so that we have more than we need in this case.

Thus we are left with the case that
$${A_1-A_1 \over A_1-A_1} \neq F_p.$$
Thus we can find $a_1,a_2,b_1,b_2$ so that for any refinement $A^{\prime} \subset A_1$
with $|A^{\prime}| \gtrsim |A_1|$, we have
$$|A_1|^2 \lesssim |(a_1-a_2) A^{\prime} + (a_1-a_2) A^{\prime} + (b_1-b_2) A^{\prime}|.$$

Now we apply Corollary 1.5, choosing $A^{\prime}$ so that
$$|(a_1-a_2) A^{\prime} + (a_1-a_2) A_1 + (b_1-b_2) A_1| \lesssim 
{|A+A| |(a_1-a_2) A_1 + (b_1-b_2) A_1| \over |A_1|}.$$
This is where we have improved over Garaev's original argument.

Then, as in the first case, estimating
$$|(a_1-a_2) A_1 + (b_1-b_2) A_1| \leq |a_1 A - a_2 A + b_1 A - b_2 A|,$$
and applying Corollary 1.4 with $X=b_0 A$ and Corollary 1.6, we obtain
$$|A_1|^3 N^4 |A|^3 \lesssim |A+A|^9.$$
Applying (2.2), we get
$$N |A|^{12} \lessapprox |A+A|^9 |AA|^3. \tag 2.5$$
Now applying (2.1), we get
$$|A|^{14} \lessapprox |A+A|^9 |AA|^4.\tag 2.6$$
Inequality (2.6) proves the Theorem. \qed \enddemo

\end
\Refs\nofrills{References}

\widestnumber\key{BGK}

\ref \key BGK \by Bourgain, J., Glibichuk, A.A., and Konyagin, S.V.
\paper Estimates for the number of sums and products and for exponential sums
in fields of prime order \jour J. London Math. Soc. (2) \vol 73 \yr 2006 \pages 380--398
\endref

\ref \key BKT \by Bourgain, J., Katz, N, and Tao, T. \paper A sum product estimate
in finite fields and Applications \jour GAFA \vol 14 \yr 2004 \pages 27--57 \endref

\ref \key G \by Garaev, M.Z. \paper An explicit sum-product estimate in $\Bbb F_p$
\jour preprint \endref

\ref \key GK \by Glibichuk, A.A, and Konyagin, S.V \paper Additive properties
of product sets in fields of prime order \jour preprint \endref

\ref \key TV \by Tao, T. and Vu, V. \book Additive Combinatorics \jour Cambridge Univ. 
Press \yr 2006 \endref







 
 \endRefs
